\documentclass[12pt, reqno]{amsart}
\usepackage{amsmath, amsthm, amscd, amsfonts, amssymb, graphicx, color}
\usepackage[bookmarksnumbered, colorlinks=false, plainpages]{hyperref}

\newtheorem{theorem}{Theorem}[section]
\newtheorem{lemma}[theorem]{Lemma}

\newtheorem{corollary}[theorem]{Corollary}
\theoremstyle{definition}
\newtheorem{definition}[theorem]{Definition}

\theoremstyle{remark}
\newtheorem{remark}[theorem]{Remark}
\numberwithin{equation}{section}

\allowdisplaybreaks
\begin{document}

\title[On $E$-frames]
{On $E$-frames}

\author[H.Hedayatirad]{Hassan Hedayatirad}
\address{Hassan Hedayatirad \\ Department of Mathematics and Computer
Sciences, Hakim Sabzevari University, Sabzevar, P.O. Box 397, IRAN}
\email{ \rm hasan.hedayatirad@hsu.ac.ir; hassan.hedayatirad67@gmail.com}
\author[T.L. Shateri]{Tayebe Lal Shateri }
\address{Tayebe Lal Shateri \\ Department of Mathematics and Computer
Sciences, Hakim Sabzevari University, Sabzevar, P.O. Box 397, IRAN}
\email{ \rm  t.shateri@hsu.ac.ir; t.shateri@gmail.com}
\thanks{*The corresponding author:
t.shateri@hsu.ac.ir; t.shateri@gmail.com (Tayebe Lal Shateri)}
 \maketitle
\begin{abstract}
In this paper, we have given some results related to this concept. In particular, we show that a classical frame can be an $E$-frame under certain conditions. Furthermore, we characterize all dual $E$-frames associated with a given $E$-frame. 
\end{abstract}
\section{Introduction}
The concept of frames for Hilbert spaces was first introduced by Duffin and Schaeffer \cite{DS}, to study non-harmonic Fourier series. It was reintroduced by Daubechies, Grossman, and Meyer \cite{DG} and has attracted mathematicians ever since. 

In this paper, we give some results on $E$-frames and dual $E$-frames. First, we recall some basic notions.

Throughout this paper, we assume that $\mathcal H$ is a separable Hilbert space.\\
Let $\mathcal X$ and $\mathcal Y$ be two sequence spaces and $E=(E_{n,k})_{n,k\geq 1}$ be an infinite matrix of real or complex numbers. We say that $E$ defines a matrix mapping from $\mathcal X$ to $\mathcal Y$, if for every sequence $x=\{x_n\}_{n=1}^{\infty}$ in $\mathcal X$, the sequence $Ex=\{(Ex)_n\}_{n=1}^{\infty}$ is in $\mathcal Y$, where
$$(Ex)_n=\sum_{k=1}^{\infty}E_{n,k}x_k,\quad n=1,2,\ldots.$$
Using the concept of matrix mapping,  the authors in \cite{TAL}, have introduced a new notion of frames which is called an $E$-frame. Let $E$ be an invertible infinite matrix mapping on $\bigoplus_{n=1}^{\infty}\mathcal{H}$. Then for each $\lbrace f_{k}\rbrace_{k=1}^{\infty}\in\bigoplus_{n=1}^{\infty}\mathcal{H}$, $$E\lbrace f_{k}\rbrace_{k=1}^{\infty}=\bigg\lbrace\sum_{k=1}^{\infty}E_{n,k}f_{k}\bigg\rbrace_{n=1}^{\infty}.$$
\begin{definition}\cite{TAL}
The sequence $\lbrace f_{k}\rbrace_{k=1}^{\infty}$ is called an $E$-frame for $\mathcal{H}$ if there exist constants $0<A\leq B<\infty$ such that 
\begin{equation}\label{$E$-frame}
A\|f\|^{2}\leq\sum_{n=1}^{\infty}\big\vert\big<f,\big(E\lbrace f_{k}\rbrace\big)_{n}\big>\big\vert^{2}\leq B\|f\|^{2}\qquad (f\in\mathcal{H}).
\end{equation}
\end{definition}
If just the right inequality in (\ref{$E$-frame}) holds, then $\lbrace f_{k}\rbrace_{k=1}^{\infty}$ is called an $E$-Bessel sequence for $\mathcal{H}$ with  $E$-Bessel bound $B$. One can easily check that for a given constant $B>0$, the sequence $\lbrace f_{k}\rbrace_{k=1}^{\infty}$ is an $E$-Bessel sequence if and only if the operator $T_{E}$ defined by $$T_{E}:\ell^{2}(\mathbb{N}):\longrightarrow\mathcal{H}\,,T\lbrace c_{k}\rbrace_{k=1}^{\infty}=\sum_{n=1}^{\infty}c_{k}\big(E\lbrace f_{k}\rbrace_{k=1}^{\infty}\big)_{n}$$ is a bounded operator from $\ell^{2}(\mathbb{N})$ to $\mathcal{H}$ with $\|T\|\leq\sqrt{B}$. We call $T_{E}$ the pre $E$-frame operator. Its adjoint, the analysis operator, is given by 
\begin{equation}\label{analysis}
T_{E}^{*}:\mathcal{H}\longrightarrow\ell^{2}(\mathbb{N})\,,T_{E}^{*}f=\bigg\lbrace\big<f,\big(E\lbrace f_{k}\rbrace_{k=1}^{\infty}\big)_{n}\big>\bigg\rbrace_{n=1}^{\infty}.
\end{equation} Composing $T_{E}$ and $T_{E}^{*}$, the $E$-frame operator $$S_{E}:\mathcal{H}\longrightarrow\mathcal{H}\,,S_{E}f=\sum_{n=1}^{\infty}\big<f,\big(E\lbrace f_{k}\rbrace_{k=1}^{\infty}\big)_{n}\big>(E\lbrace f_{k}\rbrace_{k=1}^{\infty}\big)_{n}$$ is obtained. 
\section{main results}
In this section, we assume that $\mathcal{H}$ is a separable Hilbert space and $E=(E_{n,k})_{n,k\geq 1}$ is an invertible infinite matrix mapping on $\bigoplus_{n=1}^{\infty}\mathcal{H}$. First, we show that $\bigoplus_{n=1}^{\infty}\mathcal{H}$ is a subfamily of the family of Bessel sequences in $\mathcal{H}$.

\begin{lemma}\label{lem1}
If $\lbrace f_{k}\rbrace_{k=1}^{\infty}$ belongs to $\bigoplus_{n=1}^{\infty}\mathcal{H}$, then it is a Bessel sequence.
\end{lemma}
\begin{proof}
By definition we have
\begin{equation}\label{dir}
\sum_{k=1}^{\infty}\Vert f_{k}\Vert^{2}<\infty.
\end{equation}
Now, for each $f\in\mathcal{H}$,  applying the Couchy-Shwarz inequality and \ref{dir} implies that
\begin{equation*}
\sum_{k=1}^{\infty}\left\vert\left\langle f,f_k\right\rangle\right\vert^2\leq\Vert f\Vert^2\sum_{k=1}^{\infty}\left\Vert f_k\right\Vert^2<\infty. 
\end{equation*}
This shows that $\lbrace f_k\rbrace_{k=1}^{\infty}$ is a Bessel sequence. (see \cite{ch})
\end{proof}
\begin{remark}\label{rem1}
A particular type of complex matrix is what we want to focus on now. Suppose that $E=(E_{n,k})_{n,k\geq 1}$ is an infinite complex matrix which defines an operator on $\ell^2(\mathbb{N})$ equipped with the condition that
\begin{equation}\label{Enk}
\sum_{n=1}^{\infty}\sum_{k=1}^{\infty}\vert E_{n,k}\vert^2<\infty.
\end{equation}
First note that for any sequence $\lbrace c_k\rbrace_{k=1}^{\infty}$ in $\ell^2(\mathbb{N})$,
\begin{align*}
\left\Vert E\lbrace c_k\rbrace_{k=1}^{\infty}\right\Vert_{\ell^2}^2=\sum_{n=1}^{\infty}\left\vert\sum_{k=1}^{\infty}E_{n,k}c_k\right\vert^2\leq\sum_{n=1}^{\infty}\sum_{k=1}^{\infty}\left\vert E_{n,k}\right\vert^2\sum_{k=1}^{\infty}\left\vert c_k\right\vert^2.
\end{align*}
Hence \ref{Enk} implies that $E$ is bounded.
This makes that $\lbrace E_{k,n}\rbrace_{k=1}^{\infty}$ belongs to $\ell^{2}(\mathbb{N})$ for all $n\in\mathbb{N}$. Indeed 
$$\left\lbrace E_{k,n}\right\rbrace_{k=1}^{\infty}=\left\lbrace\sum_{j=1}^{\infty}E_{k,j}\delta_{n,j}\right\rbrace_{k=1}^{\infty}=E\left\lbrace\delta_{n,k}\right\rbrace_{k=1}^{\infty},$$
where $\delta_{n,k}$ is the Kronecker delta. Therefore
\begin{equation}\label{Erow}
\left\Vert\left\lbrace E_{k,n}\right\rbrace_{k=1}^{\infty}\right\Vert_{\ell^{2}}=\left\Vert E\lbrace\delta_{n,k}\rbrace_{k=1}^{\infty}\right\Vert_{\ell^{2}}\leq\left\Vert E\right\Vert. 
\end{equation}
Applying (\ref{Erow}) to $E^{*}$ we conclude that $\lbrace E_{n,k}\rbrace_{k=1}^{\infty}$ belongs to $\ell^{2}(\mathbb{N})$ for all $n\in\mathbb{N}$.\\
Now, we will define a matrix mapping on $\bigoplus_{n=1}^{\infty}\mathcal{H}$ using $E$. For this, we assume that $\tilde{E}$ is defined as follows
\begin{equation}
\tilde{E}:\bigoplus_{n=1}^{\infty}\mathcal{H}\longrightarrow\bigoplus_{n=1}^{\infty}\mathcal{H}\quad;\quad\tilde{E}\lbrace f_k\rbrace_{k=1}^{\infty}=\left\lbrace\sum_{k=1}^{\infty}E_{n,k}f_k\right\rbrace_{n=1}^{\infty}.
\end{equation}
The above discussion together with Lemma \ref{lem1} and \cite[Corollary 3.2.5]{ch} makes $\tilde{E}$ well defined. Furthermore, for any sequence $\lbrace f_k\rbrace_{k=1}^{\infty}$ in $\bigoplus_{n=1}^{\infty}\mathcal{H}$ we see by assumption that 
\begin{equation*}
\sum_{n=1}^{\infty}\left\Vert\sum_{k=1}^{\infty}E_{n,k}f_k\right\Vert^2\leq\sum_{n=1}^{\infty}\sum_{k=1}^{\infty}\left\vert E_{n,k}\right\vert^2\sum_{k=1}^{\infty}\left\Vert f_{k}\right\Vert^2<\infty.
\end{equation*}
Thus $\tilde{E}\lbrace f_k\rbrace_{k=1}^{\infty}$ belongs to $\bigoplus_{n=1}^{\infty}\mathcal{H}$ and our aim is achieved. Therefore any infinite complex matrix which defines an operator on $\ell^2(\mathbb{N})$ and satisfies \ref{Enk} induces a matrix mapping on $\bigoplus_{n=1}^{\infty}\mathcal{H}$.
\end{remark}
Motivated by Remark \ref{rem1}, in the next result, we prove a classical frame can be an $E$-frame.
\begin{theorem}\label{th3}
Suppose $E=(E_{n,k})_{n,k\geq 1}$ is an infinite invertible complex matrix which defining an operator on $\ell^2(\mathbb{N})$ and satisfying \ref{Enk}. Suppose $\lbrace f_{k}\rbrace_{k=1}^{\infty}$ is a frame for $\mathcal{H}$ with bounds $A$ and $B$. Then $\lbrace f_{k}\rbrace_{k=1}^{\infty}$ is an $\tilde{E}$-frame with bound $\Vert E\Vert^{2}B$ and $CA$, for some $C>0$.
\end{theorem}
\begin{proof}
In the sense of Remark \ref{rem1}, let $\tilde{E}$ be the matrix mapping on $\bigoplus_{n=1}^{\infty}\mathcal{H}$ induced by $E$. We also know from Remark \ref{rem1} that $E$ is bounded. Noting that $\left\lbrace\left\langle f,f_{k}\right\rangle\right\rbrace_{k=1}^{\infty}\in\ell^{2}(\mathbb{N})$ we have
\begin{align}\label{eq1}
\sum_{n=1}^{\infty}\left\vert\left\langle f,\left(\tilde{E}\left\lbrace f_{k}\right\rbrace_{k=1}^{\infty}\right)_{n}\right\rangle\right\vert^{2}&=\sum_{n=1}^{\infty}\left\vert\left\langle f,\sum_{k=1}^{\infty}E_{n,k}f_{k}\right\rangle\right\vert^{2}=\sum_{n=1}^{\infty}\left\vert\sum_{k=1}^{\infty}\overline{E_{n,k}}\left\langle f,f_{k}\right\rangle\right\vert^{2}\nonumber\\&=\sum_{n=1}^{\infty}\left\vert\sum_{k=1}^{\infty}{E_{n,k}}\left\langle f_{k},f\right\rangle\right\vert^{2}=\left\Vert E\left\lbrace\left\langle f_{k},f\right\rangle\right\rbrace_{k=1}^{\infty}\right\Vert_{\ell^{2}}^{2}\nonumber\\&\leq\Vert E\Vert^{2}\sum_{k=1}^{\infty}\left\vert\left\langle f_{k},f\right\rangle\right\vert^{2}\leq\Vert E\Vert^{2}B\Vert f\Vert^{2}.
\end{align}
To find a lower $\tilde{E}$-frame bound for $\lbrace f_{k}\rbrace_{k=1}^{\infty}$ we note that since $E$ is invertible, so $E^{*}$ is a bijection. By \cite[Lemma 2.4.1]{ch}, there is a constant $C>0$ such that $$C\left\Vert\lbrace c_{k}\rbrace_{k=1}^{\infty}\right\Vert_{\ell^{2}}\leq\left\Vert E\lbrace c_{k}\rbrace_{k=1}^{\infty}\right\Vert_{\ell^{2}},$$ for all $\lbrace c_{k}\rbrace_{k=1}^{\infty}\in\ell^{2}(\mathbb{N})$. Consequently, the argument we stated in (\ref{eq1}) implies that
\begin{equation*}
\sum_{n=1}^{\infty}\left\vert\left\langle f,\left(\tilde{E}\left\lbrace f_{k}\right\rbrace_{k=1}^{\infty}\right)_{n}\right\rangle\right\vert^{2}=\left\Vert E\left\lbrace\left\langle f_{k},f\right\rangle\right\rbrace_{k=1}^{\infty}\right\Vert_{\ell^{2}}^{2}\geq C\sum_{k=1}^{\infty}\left\vert\left\langle f,f_{k}\right\rangle\right\vert^{2}\geq CA\Vert f\Vert^{2}.
\end{equation*}
\end{proof}
For convenience, we will use $E$ instead of $\tilde{E}$ in the notation, from now on.
\begin{corollary}
Let $E=(E_{n,k})_{n,k\geq 1}$ be an infinite diagonal matrix such that the sequence $\lbrace\lambda_n\rbrace_{n=1}^{\infty}=\lbrace E_{n,n}\rbrace_{n\geq1}$ belongs to $\ell^{\infty}\cap\ell^2(\mathbb{N})$ and $\underset{n}{inf}\vert\lambda_n\vert>0$. If $\underset{n}{sup}\vert\lambda_n\vert=\lambda$, then any frame $\lbrace f_k\rbrace_{k=1}^{\infty}$ for $\mathcal{H}$ with bounds $A$ and $B$, is an $E$-frame for $\mathcal{H}$ with bounds $CA$ and $\lambda^2B$.
\end{corollary}
\begin{proof}
First note that
\begin{equation*}
E=\begin{pmatrix}
\lambda_1&0&0&0&\cdots\\
0&\lambda_2&0&0&\cdots\\
0&0&\lambda_3&0&\cdots\\
0&0&0&\lambda_4&\cdots\\
\vdots&\vdots&\vdots&\vdots&\ddots
\end{pmatrix}
\end{equation*}
is invertible where $E^{-1}=(E_{n,k}^{-1})_{n,k\geq 1}$ is a diagonal matrix with $E_{n,n}^{-1}=\lambda_{n}^{-1}$ for each $n\in\mathbb{N}$.
Moreover for any $\lbrace c_k\rbrace_{k=1}^{\infty}\in\ell^2(\mathbb{N})$, $E\lbrace c_k\rbrace_{k=1}^{\infty}$ belongs to $\ell^2(\mathbb{N})$. Also, satisfying the conditions of \cite[Lemma 3.5.3]{ch}, $E$ defines a bounded operator on $\ell^2(\mathbb{N})$ of norm at most $\lambda$, that is for all $x\in\ell^2(\mathbb{N})$,
\begin{equation}\label{Elam}
\left\Vert Ex\right\Vert_{\ell^2}\leq\lambda\left\Vert x\right\Vert_{\ell^2}.
\end{equation}
Actually, $\lambda$ is the smallest positive number that  satisfies \ref{Elam}. In fact, suppose that $M>0$ is any arbitrary number that satisfies \ref{Elam} and take $x=\lbrace\delta_{1,k}\rbrace_{k=1}^{\infty}$ where $\delta_{1,k}$ is the Kronecker delta. Then
\begin{align*}
M=M\left\Vert\left\lbrace\delta_{1,k}\right\rbrace_{k=1}^{\infty}\right\Vert_{\ell^2}\geq\left\Vert E\left\lbrace\delta_{1,k}\right\rbrace_{k=1}^{\infty}\right\Vert_{\ell^2}=\left\Vert\left\lbrace\sum_{k=1}^{\infty}E_{n,k}\delta_{1,k}\right\rbrace_{n=1}^{\infty}\right\Vert_{\ell^2}=\left(\sum_{n=1}^{\infty}\left\vert E_{n,1}\right\vert^2\right)^{\frac{1}{2}}\leq\lambda.
\end{align*} 
This proves that $\Vert E\Vert=\lambda$. One can check easily that \ref{Enk} holds for $E$. In fact, 
\begin{align*}
\sum_{n=1}^{\infty}\sum_{k=1}^{\infty}\left\vert E_{n,k}\right\vert^2&=\sum_{n=1}^{\infty}\vert\lambda_n\vert^2<\infty.
\end{align*}
Therefore, $E$ satisfies the conditions of Theorem \ref{th3} and the proof is complete. 
\end{proof}
\begin{corollary}
Let $\lbrace f_{k}\rbrace_{k=1}^{\infty}$ be a Riesz basis in a separable Hilbert space $\mathcal{H}$ and $E=\left\lbrace\left\langle f_{k},f_{j}\right\rangle\right\rbrace_{j,k=1}^{\infty}$ be its associated Gram matrix. If $E$ satisfies \ref{Enk}, then $\lbrace f_{k}\rbrace_{k=1}^{\infty}$ is an $E$-frame for $\mathcal{H}$.
\begin{proof}
It is shown in \cite[Theorem 3.6.6]{ch} that the Gram matrix $\left\lbrace\left\langle f_{k},f_{j}\right\rangle\right\rbrace_{j,k=1}^{\infty}$ defines a bounded invertible operator on $\ell^{2}(\mathbb{N})$. The result is now completed by Theorem \ref{th3}.
\end{proof}
\end{corollary}
In the next result, we investigate, when a Bessel sequence is an $E$-Bessel sequence?
\begin{corollary}
Let $\lbrace f_{k}\rbrace_{k=1}^{\infty}$ be a Bessel sequence on $\mathcal{H}$. Suppose that $E=(E_{n,k})_{n,k\geq 1}$ is an infinite complex matrix satisfying the conditions of Theorem \ref{th3}. Then $\lbrace f_{k}\rbrace_{k=1}^{\infty}$ is an $E$-Bessel sequence if and only if $\overline{E}(T^{*}f)\in\ell^{2}(\mathbb{N})$, where $T^{*}$ is the analysis operator of $\lbrace f_{k}\rbrace_{k=1}^{\infty}$.
\end{corollary}
\begin{proof}
Note that $T^{*}$ is well defined since $\lbrace f_{k}\rbrace_{k=1}^{\infty}$ is a Bessel sequence. Now for given $f\in\mathcal{H}$ we have
\begin{align*}
\sum_{n=1}^{\infty}\left\Vert\left\langle f,\left(E\left\lbrace f_{k}\right\rbrace_{k=1}^{\infty}\right)_{n}\right\rangle\right\Vert^{2}&=\sum_{n=1}^{\infty}\left\Vert\left\langle f,\sum_{k=1}^{\infty}E_{n,k}f_{k}\right\rangle\right\Vert^{2}=\sum_{n=1}^{\infty}\left\Vert\sum_{k=1}^{\infty}\overline{E}_{n,k}\left\langle f,f_{k}\right\rangle\right\Vert^{2}\\&=\sum_{n=1}^{\infty}\left\Vert\left(\overline{E}\left\lbrace\left\langle f,f_{k}\right\rangle\right\rbrace_{k=1}^{\infty}\right)_{n}\right\Vert^{2}=\left\Vert\overline{E}\left\lbrace\left\langle f,f_{k}\right\rangle\right\rbrace_{k=1}^{\infty}\right\Vert_{\ell^2}^{2}.
\end{align*}
\end{proof}
Similar to \cite[Proposition 5.5.8]{ch}, we can prove the following theorem.
\begin{theorem}
Let $\lbrace f_{k}\rbrace_{k=1}^{\infty}$ be a frame with bounds $A$,$B$. Suppose that the infinite complex matrix $E=\lbrace E_{n,k}\rbrace_{k,n\in\mathbb{N}}$ defines a matrix mapping on $\bigoplus_{n=1}^{\infty}\mathcal{H}$ and satisfies the two conditions
\begin{align*}
&b:=\underset{k}{sup}\sum_{j=1}^{\infty}\left\vert\sum_{n=1}^{\infty}E_{n,k}\overline{E}_{n,j}\right\vert<\infty,\\&a:=\underset{k}{inf}\left(\sum_{n=1}^{\infty}\left\vert E_{n,k}\right\vert^{2}-\sum_{j\neq K}\left\vert\sum_{n=1}^{\infty}E_{n,k}\overline{E}_{n,j}\right\vert\right)>0.
\end{align*}
Then $\lbrace f_{k}\rbrace_{k=1}^{\infty}$ is an $E$-frame with bounds $aA$ and $bB$.
\end{theorem}

In the following, we give the definition of an $E$-basis from \cite{TAL}.
\begin{definition}
The sequence $\lbrace g_{k}\rbrace_{k=1}^{\infty}$ is an $E$-basis for $\mathcal{H}$ if for each $f\in\mathcal{H}$ there exists a unique scalar coefficients $\lbrace c_{k}(f)\rbrace_{k=1}^{\infty}$ such that 
\begin{equation}\label{ebas}
f=\sum_{k=1}^{\infty}c_{k}\left(f\right)\left(E\left\lbrace g_{j}\right\rbrace_{j=1}^{\infty}\right)_{k}.
\end{equation}
Also, an $E$-basis $\lbrace g_{k}\rbrace_{k=1}^{\infty}$ is $E$-orthonormal if $\lbrace g_{k}\rbrace_{k=1}^{\infty}$ is an $E$-orthonormal system. i.e., if 
\end{definition}
\begin{equation*}
\left\langle\left(E\left\lbrace g_{j}\right\rbrace_{j=1}^{\infty}\right)_{n},\left(E\left\lbrace g_{j}\right\rbrace_{j=1}^{\infty}\right)_{k}\right\rangle=\delta_{n,k}=\begin{cases}
1&n=k,\\0&n\neq k.
\end{cases}
\end{equation*}
Suppose that $\lbrace e_{k}\rbrace_{k=1}^{\infty}$ is an orthonormal basis for $\mathcal{H}$. Then  it is easy to show that for every invertible infinite matrix operator $E$ on $\bigoplus_{n=1}^{\infty}\mathcal{H}$, $E^{-1}\lbrace e_{k}\rbrace_{k=1}^{\infty}$ is an $E$-orthonormal basis for $\mathcal{H}$. So every separable Hilbert space has an $E$-orthonormal basis. Furthermore, every $E$-orthonormal basis is an $E$-Bessel sequence \cite{TAL}. If $\lbrace g_{k}\rbrace_{k=1}^{\infty}$ is an $E$-orthonormal basis for $\mathcal{H}$, then for every $f\in\mathcal{H}$ there exist unique coefficients $\lbrace c_{n}(f)\rbrace_{k=1}^{\infty}$ in $\ell^{2}(\mathbb{N})$ such that 
\begin{equation}
f=\sum_{n=1}^{\infty}c_{n}\left(f\right)\left(E\left\lbrace g_{k}\right\rbrace_{k=1}^{\infty}\right)_{n}.
\end{equation}
In fact $$c_{m}(f)=\left\langle f,\left(E\left\lbrace g_{k}\right\rbrace_{k=1}^{\infty}\right)_{m}\right\rangle\qquad(m\in\mathbb{N}).$$
By a result proved by Cassaza, in a complex Hilbert space, every frame can be written as a multiple of three orthonormal bases. We refer to \cite{ch} for the statement and the proof. Here, We have given an interpretation of this result in terms of $E$-frames and $E$-orthonormal bases.
\begin{theorem}
Suppose $\mathcal{H}$ is a complex Hilbert space and that $\lbrace f_{k}\rbrace_{k=1}^{\infty}$ is an $E$-frame for $\mathcal{H}$ with synthesis operator $T_{E}$. Then for every $\varepsilon\in(0,1)$ there exists three $E$-orthonormal bases $\lbrace g_{k}^{1}\rbrace_{k=1}^{\infty}$, $\lbrace g_{k}^{2}\rbrace_{k=1}^{\infty}$ and $\lbrace g_{k}^{3}\rbrace_{k=1}^{\infty}$ such that 
\begin{equation*}
E\lbrace f_{k}\rbrace_{k=1}^{\infty}=\dfrac{\Vert T\Vert}{1-\varepsilon}\left(E\lbrace g_{k}^{1}\rbrace_{k=1}^{\infty}+E\lbrace g_{k}^{2}\rbrace_{k=1}^{\infty}+E\lbrace g_{k}^{3}\rbrace_{k=1}^{\infty}\right),
\end{equation*}
where $T$ is a bounded linear surjection on $\mathcal{H}$.
\end{theorem}
\begin{proof}
Let $\lbrace g_{k}\rbrace_{k=1}^{\infty}$ be an $E$-orthonormal basis for $\mathcal{H}$.  Suppose that $\lbrace \delta_{n}\rbrace_{n=1}^{\infty}$ be the canonical orthonormal basis for $\ell^{2}(\mathbb{N})$. Consider the mapping 
\begin{equation*}
\varphi:\mathcal{H}\longrightarrow\ell^{2}(\mathbb{N})\,;\,\varphi\left(f\right)=\left\lbrace \left\langle f,\left(E\left\lbrace g_{k}\right\rbrace_{k=1}^{\infty}\right)_{n}\right\rangle\right\rbrace_{n=1}^{\infty}
\end{equation*}
and set $T=T_{E}\varphi$. It is easy to show that $\varphi$ is an isometric isomorphism. Thus $T$ defines a bounded linear surjection on $\mathcal{H}$. Given $\varepsilon\in(0,1)$, consider the operator 
\begin{equation}\label{eqD}
D:\mathcal{H}\longrightarrow\mathcal{H}\,;\,D:=\dfrac{1}{2}id_{\mathcal{H}}+\dfrac{(1-\varepsilon)T}{2\Vert T\Vert}.
\end{equation}
Clearly $\Vert id_{\mathcal{H}}-D\Vert<1$ and so $D$ is invertible. Hence \cite[Lemma 2.4.6]{ch} implies that $D=\dfrac{V(W+W^{*})}{2}$, where $V$ and $W$ are unitary operators. Now, by (\ref{eqD}) we can now write
\begin{equation*}
T=\dfrac{2\Vert T\Vert(D-id_{\mathcal{H}})}{1-\varepsilon}=\dfrac{\Vert T\Vert}{1-\varepsilon}\left(VW+VW^{*}-id_{\mathcal{H}}\right).
\end{equation*}
It follows from here that for each $n\in\mathbb{N}$,
\begin{align*}
\left(E\left\lbrace f_{k}\right\rbrace_{k=1}^{\infty}\right)_{n}=T_{E}\delta_{n}&=T\left(E\left\lbrace g_{k}\right\rbrace_{k=1}^{\infty}\right)_{n}\\&=\dfrac{\Vert T\Vert}{1-\varepsilon}\left(VW\left(E\left\lbrace g_{k}\right\rbrace_{k=1}^{\infty}\right)_{n}+VW^{*}\left(E\left\lbrace g_{k}\right\rbrace_{k=1}^{\infty}\right)_{n}-\left(E\left\lbrace g_{k}\right\rbrace_{k=1}^{\infty}\right)_{n}\right)\\&=\dfrac{\Vert T\Vert}{1-\varepsilon}\left(\left(E\left\lbrace VWg_{k}\right\rbrace_{k=1}^{\infty}\right)_{n}+\left(E\left\lbrace VW^{*}g_{k}\right\rbrace_{k=1}^{\infty}\right)_{n}-\left(E\left\lbrace g_{k}\right\rbrace_{k=1}^{\infty}\right)_{n}\right).
\end{align*}
Since $VW$ and $VW^{*}$ are unitary operators, the result now follows from \cite[Theorem 2.4]{TAL}.
\end{proof}


\end{document}